\documentclass[11pt]{article}

\usepackage{amssymb}

\def\G{{\cal G}}
\def\C{{\mathbb C}}
\def\N{{\mathbb N}}
\def\Z{{\mathbb Z}}
\def\R{{\mathbb R}}

\def\T{{\mathbb T}}

\def\Mat{{\rm Mat}}
\def\rank{{\rm rank}\,}
\def\CPA{{CPA}}
\def\eps{\varepsilon}
\def\Ad{{\rm Ad}\,}
\def\id{{\rm id}}
\def\tr{{\rm tr}}

\def\Aut{\mathop{\mbox{\rm Aut\,}}}

\def\rcp{{rcp}}

\makeatletter
\@addtoreset{equation}{section}

\renewcommand{\subsection}{\noindent\refstepcounter{subsection}
\thesubsection.\ }

\newtheorem{theorem}{Theorem}[section]
\newtheorem{lemma}[theorem]{Lemma}

\newtheorem{proposition}[theorem]{Proposition}
\newtheorem{example}[theorem]{Example}

\newenvironment{verif}[1]{\medskip\noindent{\bf Proof#1.}}
    {\newline\mbox{\ }\hfill\rule{0.5em}{0.5em}\smallskip}

\begin{document}

\title{Entropy of automorphisms of ${\rm II}_1$-factors
arising from the dynamical systems theory}

\author{V.Ya. Golodets, S.V. Neshveyev}

\date{\it B. Verkin Institute for Low Temperature Physics and
Engineering, National Academy of Sciences of Ukraine,
47, Lenin Ave., 310164, Kharkov, Ukraine}

\maketitle

\begin{abstract}
Let a countable amenable group $G$ acts freely and ergodically on
a Lebesgue space $(X, \mu)$, preserving the measure
$\mu$. If $T\in\Aut(X, \mu)$ is an automorphism of the equivalence
relation defined by $G$ then $T$ can be extended to an
automorphism $\alpha_T$ of the II$_1$-factor
$M=L^\infty(X,\mu)\rtimes G$. We prove that if $T$ commutes with
the action of $G$ then $H(\alpha_T)=h(T)$, where $H(\alpha_T)$ is
the Connes-St{\o}rmer entropy of $\alpha_T$, and
$h(T)$ is the Kolmogorov--Sinai entropy of $T$. We prove also that
for given $s$ and $t$, $0\le s\le t\le\infty$, there exists a
$T$ such that $h(T)=s$ and $H(\alpha_T)=t$.
\end{abstract}

\bigskip

\section*{Introduction}

Entropy is an important notion in classical statistical mechanics
and information theory. Initially the conception of entropy for
automorphism in the ergodic theory was introduced by Kolmogorov
and Sinai in 1958. This invariant proved to be extremely useful in
the classical dynamical systems theory and topological dynamics.
The extension of this notion onto quantum dynamical systems was
done by Connes, Narnhofer, St{\o}rmer and Thirring~\cite{CS,CNT}.
At the present time there are several other promising approaches
to entropy of $C^*$-dynamical systems~\cite{S, AF, V}.

An important trend in dynamical entropy is its computation for
various models. A lot of interesting results was obtained in this
field in the recent years. We note several of them. St{\o}rmer,
Voiculescu~\cite{SV}, and the second author~\cite{N} computed the
entropy of Bogoliubov automorphisms of CAR and CCR algebras (see
also~\cite{BG,GN2}). Pimsner, Popa~\cite{PP}, Choda~\cite{Ch1}
computed the entropy of shifts of Temperley-Lieb algebras,
Choda~\cite{Ch2}, Hiai~\cite{H} and St{\o}rmer~\cite{St} computed
the entropy of canonical shifts. The first author,
St{\o}rmer~\cite{GS1,GS2}, Price~\cite{Pr} computed entropy for a
wide class of binary shifts.

In this paper we consider automorphisms of II$_1$ factors arising
from the dynamical systems theory. Let a countable group
$G$ acts freely and ergodically on a Lebesgue space $(X, \mu)$ and
preserves $\mu$. Then one can construct the crossed product
$M=L^\infty(X, \mu)\rtimes G$, which, as well known, is a II$_1$-factor. If
$T\in\Aut(X, \mu)$ defines an automorphism of the ergodic
equivalence relation induced by $G$ then $T$ can be extended to an
automorphism $\alpha_T$ of $M$~\cite{FM}. It is a natural problem
to compute the dynamical entropy $H(\alpha_T)$ in the sense
of~\cite{CS} and to compare it with the Kolmogorov-Sinai entropy
$h(T)$ of $T$. It should be noted that this last problem is a part of a more
general problem. Namely, let $M$ be a II$_1$-factor,
$\alpha\in\Aut M$,
$A$ its $\alpha$-invariant Cartan subalgebra, $\alpha(A)=A$, then
it is nature to investigate when $H(\alpha)$ is equal to
$H(\alpha|_A)$.

These problems are studied in our paper. In Section~\ref{2} we
prove that if $T$ commutes with the action of $G$ then
$H(\alpha_T)=h(T)$. More generally, we prove that this result is
valid for crossed products of arbitrary algebras for entropies of
Voiculescu~\cite{V} and of Connes-Narnhofer-Thirring~\cite{CNT}.
In Section~\ref{3} we consider two examples to illustrate this
result. These examples give non-isomorphic ergodic automorphisms
of the hyperfinite ergodic equivalence relation with the same
entropy. In Section~\ref{4} we construct several examples showing
that the entropies $h(T)$ and $H(\alpha_T)$ can be distinct. These
systems are non-commutative analogs of dynamical systems of
algebraic origin (see \cite{A,Y,LSW,S}). In particular, some of our
examples are automorphisms of non-commutative tori.
In Section~\ref{5} we construct flows $T_t$
such that $H(\alpha_{T_1})>h(T_1)$. In particular, we show that
the values $h(T)$ and $H(\alpha_T)$ can be arbitrary.

\bigskip\bigskip

\section{Computation of entropy of automorphisms of crossed
products}\label{2}


Let $(X, \mu)$ be a Lebesgue space, $G$ a countable amenable
group of automorphisms $S_g$, $g\in G$, of $(X,\mu)$ preserving
$\mu$, and $T$ an automorphism of $(X, \mu)$, $\mu\circ T=\mu$, such
that
$$
TS_g=S_gT,\quad g\in G.
$$

\begin{theorem} \label{2.1}
Let $(X, \mu)$, $G$ and $T$ be as above. Suppose $G$ acts freely
and ergodically on $(X, \mu)$. Then $M=L^\infty(X, \mu)\rtimes_SG$
is the hyperfinite II$_1$-factor with the trace-state $\tau$
induced by $\mu$. The automorphism $T$ can be canonically extended
to an automorphism $\alpha_T$ of $M$, and $$ H(\alpha_T)=h(T)\,,
$$ where $H(\alpha_T)$ is the Connes-St{\o}rmer entropy of
$\alpha_T$, and $h(T)$ is the Kolmogorov-Sinai entropy of $T$.
\end{theorem}

We will prove the following more general result.

\begin{theorem} \label{2.2}
Let $M$ be an approximately finite-dimensional W$^*$-algebra,
$\sigma$ its normal state, $T$ a $\sigma$-preserving automorphism.
Suppose a discrete amenable group $G$ acts on $M$ by automorphisms
$S_g$ that commute with $T$ and preserve $\sigma$. The
automorphism $T$ defines an automorphism $\alpha_T$ of
$M\rtimes_SG$, and the state $\sigma$ is extended to the dual
state which we continue to denote by $\sigma$. Then

(i) $hcpa_\sigma(\alpha_T)=hcpa_\sigma(T)$,
where $hcpa_\sigma$ is the completely positive approximation entropy of
Voiculescu \cite{V};

(ii) $h_\sigma(\alpha_T)=h_\sigma(T)$, where $h_\sigma$ is the
dynamical entropy of Connes-Narnhofer-Thirring~\cite{CNT}.
\end{theorem}

Since CNT-entropy coincides with KS-entropy in the classical case, and with
CS-entropy for tracial $\sigma$ and approximately finite-dimensional $M$,
Theorem~\ref{2.1} follows from Theorem~\ref{2.2}.

\smallskip

To prove Theorem~\ref{2.2} we will generalize a construction
of Voiculescu~\cite{V}.

\begin{lemma} \label{2.3}
Let $B$ be a C$^*$-algebra, $x_1,\ldots,x_n\in B$. Then the mapping
$\Psi\colon\Mat_n(\C)\otimes B\to B$,
$$
\Psi(e_{ij}\otimes b)=x_ibx^*_j,
$$
is completely positive.
\end{lemma}

\begin{verif}{}
Consider the element $V\in\Mat_n(B)=\Mat_n(\C)\otimes B$,
$$
V=\pmatrix{x_1 & \ldots & x_n\cr
       0   & \ldots & 0  \cr
           & \ldots &    \cr
       0   & \ldots & 0  \cr}.
$$
Consider also the projection $p=e_{11}\otimes 1\in\Mat_n(\C)\otimes B$.
Then $\Psi$ is the mapping
$\Mat_n(B)\to p\Mat_n(B)p=B$, $x\mapsto VxV^*$.
\end{verif}

Let $\lambda$ be the canonical representation of $G$ in $M\rtimes
G$, so that $(\Ad\lambda(g))(a)=S_g(a)$ for $a\in M$.

\begin{lemma} \label{2.4}
For any finite subset $F$ of $G$,
there exist normal unital completely positive mappings
$I_F\colon B(l^2(F))\otimes M\to M\rtimes G$ and
$J_F\colon M\rtimes G\to B(l^2(F))\otimes M$ such that
\begin{eqnarray*}
I_F(e_{g,h}\otimes a)
 &=&{1\over|F|}\lambda(g)a\lambda(h)^*
     ={1\over|F|}\lambda(gh^{-1})S_h(a),\\
J_F(\lambda(g)a)
 &=&\sum_{h\in F\cap g^{-1}F}e_{gh,h}\otimes S_{h^{-1}}(a),\\
(I_F\circ J_F)(\lambda(g)a)
 &=&{|F\cap g^{-1}F|\over|F|}\lambda(g)a,\\
\sigma\circ I_F
 &=&\tr_F\otimes\sigma,\ \ \alpha_T\circ I_F=I_F\circ(\id\otimes T),\\
(\tr_F\otimes\sigma)\circ J_F
 &=&\sigma,\ \ (\id\otimes T)\circ J_F=J_F\circ\alpha_T,
\end{eqnarray*}
where $\tr_F$ is the unique tracial state on $B(l^2(F))$.
\end{lemma}

\begin{verif}{}
The complete positivity of $I_F$ follows from Lemma~\ref{2.3}.
Consider $J_F$. Suppose that $M\subset B(H)$, and consider
the regular representation of $M\rtimes G$ on $l^2(G)\otimes H$:
$$
\lambda(g)(\delta_h\otimes\xi)=\delta_{gh}\otimes\xi, \ \
 a(\delta_h\otimes\xi)=\delta_h\otimes S_{h^{-1}}(a)\xi\ \ (a\in M).
$$
Let $P_F$ be the projection onto $l^2(F)\otimes H$. Then
a direct computation shows that the mapping $J_F(x)=P_FxP_F$,
$x\in M\rtimes G$, has the form written above. All others assertions follow
immediately.
\end{verif}

\begin{verif}{ of Theorem~\ref{2.2}}

\noindent (i) Since there exists a $\tau$-preserving conditional expectation
$M\rtimes G\to M$, we have $hcpa_\sigma(\alpha_T)\ge hcpa_\sigma(T)$.
To prove the opposite inequality we have to show that
$hcpa_\sigma(\alpha_T,\omega)\le hcpa_\sigma(T)$ for any finite subset
$\omega$ of $M\rtimes G$. Fix $\eps>0$. We can find a finite subset $F$ of
$G$ such that $||(I_F\circ J_F)(x)-x||_\sigma<\eps$ for any $x\in\omega$.
Let $(\psi,\phi,B)\in\CPA(B(l^2(F))\otimes M,\tr_F\otimes\sigma)$.
Then $(I_F\circ\psi,\phi\circ J_F,B)\in\CPA(M\rtimes G,\sigma)$.
Suppose
$$
||(\psi\circ\phi)(J_F(x))-J_F(x)||_{\tr_F\otimes\sigma}<\delta
$$
for some $x\in\alpha^k_T(\omega)$ and $k\in\N$. Then
$$
||(I_F\circ\psi\circ\phi\circ J_F)(x)-x||_\sigma
 \le||(\psi\circ\phi)(J_F(x))-J_F(x)||_{\sigma\circ I_F}
      +||(I_F\circ J_F)(x)-x||_\sigma<\delta+\eps,
$$ where we have used the facts that $\sigma\circ
I_F=\tr_F\otimes\sigma$ and that $\alpha_T$ commutes with
$I_F\circ J_F$. Since $J_F\circ\alpha_T=(\id\otimes T)\circ J_F$,
we infer that $$
\rcp_\sigma(\omega\cup\alpha_T(\omega)\cup\ldots\cup\alpha^{n-1}_T(\omega);
   \delta+\eps)\le\rcp_{\tr_F\otimes\sigma}(J_F(\omega)\cup
      \ldots\cup(\id\otimes T)^{n-1}(J_F(\omega));\delta),
$$
so that (for $\delta<\eps$)
\begin{eqnarray*}
hcpa_\sigma(\alpha_T,\omega;2\eps)
 &\le& hcpa_\sigma(\alpha_T,\omega;\eps+\delta)
        \le hcpa_{\tr_F\otimes\sigma}(\id\otimes T,J_F(\omega);\delta)\\
 &\le& hcpa_{\tr_F\otimes\sigma}(\id\otimes T)=hcpa_\sigma(T),
\end{eqnarray*}
where the last equality follows from the subadditivity of the
entropy~\cite{V}. Since $\eps>0$ was arbitrary, the proof of the
inequality $hcpa_\sigma(\alpha_T,\omega)\le hcpa_\sigma(T)$ is
complete.

\smallskip

\noindent (ii) We always have $h_\sigma(\alpha_T)\ge h_\sigma(T)$.
To prove the opposite inequality consider a channel $\gamma\colon
B\to M\rtimes G$, i.~e., a unital completely positive mapping of a
finite-dimensional C$^*$-algebra~$B$. We have to prove that
$h_\sigma(\alpha_T;\gamma)\le h_\sigma(T)$. Fix $\eps>0$. We can
choose $F$ such that
$$
||(I_F\circ J_F\circ\gamma-\gamma)(x)||_\sigma\le\eps||x||
 \ \ \hbox{for any}\ x\in B.
$$
By~\cite[Theorem IV.3]{CNT},
\begin{equation} \label{e2.1}
{1\over n}H_\sigma(\gamma,\alpha_T\circ\gamma,\ldots,\alpha^{n-1}_T\circ\gamma)
   \le\delta+{1\over n}
       H_\sigma(I_F\circ J_F\circ\gamma,\alpha_T\circ I_F\circ J_F\circ\gamma,
          \ldots,\alpha^{n-1}_T\circ I_F\circ J_F\circ\gamma),
\end{equation}
where $\delta=\delta(\eps,\rank B)\to0$ as $\eps\to0$. Since
$\sigma\circ I_F=\tr_F\otimes\sigma$, it is easy to see from
the definition of mutual entropy $H_\sigma$ \cite{CNT} that
\begin{equation} \label{e2.2}
H_\sigma(I_F\circ J_F\circ\gamma,I_F\circ J_F\circ\alpha_T\circ\gamma,
   \ldots,I_F\circ J_F\circ\alpha^{n-1}_T\circ\gamma)
 \le H_{\tr_F\otimes\sigma}(J_F\circ\gamma,J_F\circ\alpha_T\circ\gamma,
        \ldots,J_F\circ\alpha^{n-1}_T\circ\gamma)
\end{equation}
Since $I_F\circ J_F$ commutes with $\alpha_T$, and
$J_F\circ\alpha_T=(\id\otimes T)\circ J_F$, we infer from (\ref{e2.1}) and
(\ref{e2.2}) that
$$
h_\sigma(\alpha_T;\gamma)\le\delta
 +h_{\tr_F\otimes\sigma}(\id\otimes T;J_F\circ\gamma)
 \le\delta+h_{\tr_F\otimes\sigma}(\id\otimes T).
$$
Since we could choose $F$ such that $\delta$ was
arbitrary small, we see that it suffices to prove that
$h_{\tr_F\otimes\sigma}(\id\otimes T)=h_\sigma(T)$. For abelian $M$
this is proved by standard arguments, using
\cite[Corollary VIII.8]{CNT}. To handle the general case we need the following
lemma.

\begin{lemma} \label{2.5}
For any finite-dimensional C$^*$-algebra $B$, any state $\phi$ of $B$, and
any positive linear functional $\psi$ on $\Mat_n(\C)\otimes B$, we have
$$
S(\tr_n\otimes\phi,\psi)\le S(\phi,\psi|_B)+2\psi(1)\log n.
$$
\end{lemma}

\begin{verif}{}
By \cite[Theorem 1.13]{OP},
$$
S(\tr_n\otimes\phi,\psi)=S(\phi,\psi|_B)+S(\psi\circ E,\psi),
$$
where $E=\tr_n\otimes\id\colon\Mat_n(\C)\otimes B\to B$ is the
($\tr_n\otimes\phi$)-preserving conditional expectation (note that we adopt
the notations of~\cite{CNT}, so we denote by $S(\omega_1,\omega_2)$ the
quantity which is denoted by $S(\omega_2,\omega_1)$ in~\cite{OP}). By the
Pimsner-Popa inequality~\cite[Theorem 2.2]{PP}, we have
$$
E(x)\ge{1\over n^2}x\ \ \hbox{for any}\ \ x\in\Mat_n(\C)\otimes B,\ x\ge0.
$$
In particular, $\psi\circ E\ge{1\over n^2}\psi$, whence
$S(\psi\circ E,\psi)\le2\psi(1)\log n$.
\end{verif}

Since $M$ is an AFD-algebra, to compute the entropy of $\id\otimes T$ it
suffices to consider subalgebras of the form $B(l^2(F))\otimes B$, where
$B\subset M$. From Lemma~\ref{2.5} and the definitions~\cite{CNT} we
immediately get
$$
h_{\tr_F\otimes\sigma}(\id\otimes T;B(l^2(F))\otimes B)\le h_\sigma(T;B)
 +2\log|F|.
$$
Hence
$h_{\tr_F\otimes\sigma}(\id\otimes T)\le h_\sigma(T)+2\log|F|$.
Applying this inequality to $T^m$, we obtain
$$
h_{\tr_F\otimes\sigma}((\id\otimes T)^m)\le h_\sigma(T^m)+2\log|F|
   \ \ \forall m\in\N.
$$
But since $M$ is an AFD-algebra, we have
$h_{\tr_F\otimes\sigma}((\id\otimes T)^m)
   =m\cdot h_{\tr_F\otimes\sigma}(\id\otimes T)$ and
$h_\sigma(T^m)=m\cdot h_\sigma(T)$. So dividing the above inequality by $m$,
and letting $m\to\infty$, we obtain
$h_{\tr_F\otimes\sigma}(\id\otimes T)\le h_\sigma(T)$, and the proof of
Theorem is complete.
\end{verif}

\noindent{\bf Remarks.}

\noindent(i) For any AFD-algebra $N$ and any normal state $\omega$ of $N$,
we have $h_{\omega\otimes\sigma}(\id\otimes T)=h_\sigma(T)$. Indeed, we may
suppose that $N$ is finite-dimensional and $\omega$ is faithful (because
if $p$ is the support of $\omega$, then $h_{\omega\otimes\sigma}(\id\otimes T)
=h_{\omega\otimes\sigma}((\id\otimes T)|_{pNp\otimes M})$). Now the only
thing we need is a generalization of the Pimsner-Popa inequality. Let
$p_1,\ldots,p_m$ be the atoms of a maximal abelian subalgebra of the
centralizer of the state $\omega$. Then
$$
(\omega\otimes\id)(x)\ge\left(\sum^m_{i=1}{1\over\omega(p_i)}\right)^{-1}x
   \ \ \hbox{for any}\ \ x\in N\otimes M,\ x\ge0,
$$
by \cite[Theorem 4.1 and Proposition 5.4]{L}.

\noindent(ii) By Corollary 3.8 in \cite{V}, $hcpa_\mu(T)=h(T)$ for ergodic
$T$. For non-ergodic $T$, the entropies can be distinct. Indeed, let $X_1$
be a $T$-invariant measurable subset of $X$, $\lambda=\mu(X)$, $0<\lambda<1$.
Set $\mu_1=\lambda^{-1}\mu|_{X_1}$, $T_1=T|_{X_1}$, $X_2=X\backslash X_1$,
$\mu_2=(1-\lambda)^{-1}\mu|_{X_2}$, $T_2=T|_{X_2}$. It is easy to see that
$h(T)=\lambda h(T_1)+(1-\lambda)h(T_2)$. On the other hand, it can be proved
that
$$
hcpa_\mu(T)=\max\{hcpa_{\mu_1}(T_1),hcpa_{\mu_2}(T_2)\}.
$$
So if $h(T_1),h(T_2)<\infty$, $h(T_1)\ne h(T_2)$, then $h(T)<hcpa_\mu(T)$.

To obtain an invariant which coincides with KS-entropy in the classical case,
one can modify Voiculescu's definition replacing $\rank B$ with
$\exp S(\sigma\circ\psi)$ in~\cite[Definition 3.1]{V}. Theorem~\ref{2.2}
remains true for this modified entropy.

\bigskip\bigskip

\section{Examples}\label{3}

We present two examples to illustrate Theorem~\ref{2.1}. These examples
give non-isomorphic ergodic automorphisms of amenable equivalence
relations with the same KS-entropy.\medskip

Let us first describe a general construction.

\begin{proposition} \label{3.1}
Let $S_0$, $S_1$, $S_2$ be ergodic automorphisms of $(X,\mu)$ such
that $S_0$ commutes with $S_1$ and $S_2$, and $S_1$ is conjugate
with neither $S_2$, nor $S^{-1}_2$ by an automorphism commuting
with $S_0$. Set $M_i=L^\infty(X,\mu)\rtimes_{S_i}\Z$, $i=1,2$, and
let $\alpha_i$ be the automorphism of $M_i$ induced by $S_0$. Then
there is no isomorphism $\phi$ of $M_1$ onto $M_2$ such that
$\phi\circ\alpha_1=\alpha_2\circ\phi$ and
$\phi(L^\infty(X,\mu))=L^\infty(X,\mu)$.
\end{proposition}

\begin{verif}{}
Suppose such a $\phi$ exists. Let $U_i\in M_i$ be the unitary corresponding
to $\alpha_i$, $i=1,2$, $A=L^\infty(X,\mu)\subset M_1$.
Set $U=\phi^{-1}(U_2)$.
Since $U$ is a unitary operator from $M_1$ such that $(\Ad U)(A)=A$,
it is easy to check that $U$ has the form
$$
U=\sum_{i\in\Z}a_iU_1^i E_i,\quad a_i\in\T,
$$
where $\{E_i\} $ is a family of projections from $A$, $E_iE_j=0$, for
$i\not=j$, $\sum\limits_iE_i=\sum\limits_iU_1^iE_iU_1^{-i}=I$.
Since $\alpha_1(U)=U$, we have $\alpha_1(E_i)=E_i$, $i\in\Z$.
But $S_0$ is ergodic, therefore $E_i=I$ or $E_i=0$. Hence
$U=a_iU_1^i$ for some $i\in\Z$ and $a_i\in\T$. Since $\phi$ is an isomorphism,
we have either $i=-1$, or $i=1$. We see that $\phi|_{L^\infty(X,\mu)}$ is an
automorphism that commutes with $S_0$ and conjugates $S_2$ with either
$S^{-1}_1$, or $S_1$.
\end{verif}

\noindent{\bf Remark.} It follows from Proposition~\ref{3.1} that $S_0$
defines non-isomorphic automorphisms of the ergodic equivalence
relations induced by $S_1$ and $S_2$ on $X$ correspondingly, despite
of $H(\alpha_1)=H(\alpha_2)=h(S_0)$.

\begin{example} \label{3.2}
\rm
Let $X=[0, 1]$ be the unit interval,
$\mu$ the Lebesgue measure on $X$, $t_0$, $t_1$ and $t_2$ irrational numbers
from $[0, 1]$ such that $t_2\ne t_1,\,1-t_1$.
Consider the shifts $S_ix=x+t_i\pmod1$, $x\in[0, 1]$. Any automorphism
of $X$ commuting with $S_0$ commutes with $S_1$ and $S_2$. Since
$S_1\ne S^{\pm1}_2$, Proposition~\ref{3.1} is applicable. Note that
$h(S_0)=0$.
\end{example}

\begin{example} \label{3.3}
\rm Let $(X, \mu)$ be a Lebesgue space, $T_t$ a Bernoulli flow on
$(X, \mu)$ with $h(T_1)=\log2$ \cite{O}. Choose $t_i\in\R$,
$t_i\ne0$ ($i=0,1,2$), $t_1\ne\pm t_2$, and set $S_i=T_{t_i}$.
Then $h(S_1)\ne h(S_2)$, and we can apply Proposition~\ref{3.1}.
\end{example}

\bigskip\bigskip

\section{Entropy of automorphisms and their restrictions to a
Cartan subalgebra} \label{4}

Let $M$ be a II$_1$-factor, $A$ its Cartan subalgebra,
$\alpha\in\Aut M$ such that $\alpha(A)=A$. We consider cases when
$H(\alpha)>H(\alpha|_A)$.

\medskip

Suppose a discrete abelian group $G$ acts freely and ergodically by
automorphisms $S_g$ on $(X,\mu)$, $\beta$ an automorphism of $G$, and
$S$ an automorphism of $(X,\mu)$ such that $TS_g=S_{\beta(g)}T$. Then
$T$ induces an automorphism $\alpha_T$ of $M=L^\infty(X,\mu)\rtimes_SG$.
Explicitly,
$$
\alpha_T(f)(x)=f(T^{-1}x) \ \hbox{for}\ f\in L^\infty(X,\mu),
 \ \alpha_T(\lambda(g))=\lambda(\beta(g)).
$$
The algebra $A=L^\infty(X,\mu)$ is a Cartan subalgebra of $M$. On the other
hand, the operators $\lambda(g)$ generate a maximal abelian subalgebra
$B\cong L^\infty(\hat G)$ of $M$, and $\alpha_T|_B=\hat\beta$, the dual
automorphism of $\hat G$. We have
$$
H(\alpha_T)\ge\max\{h(T),h(\hat\beta)\},
$$
so if $h(\hat\beta)>h(T)$, then $H(\alpha_T)>H(\alpha_T|_A)$.

\medskip

To construct such examples we consider systems of algebraic origin.

Let $G_1$ and $G_2$ be discrete abelian groups, and $T_1$ an
automorphism of $G_1$. Suppose there exists an embedding $l\colon
G_2\hookrightarrow\hat G_1$ such that $l(G_2)$ is a dense $\hat
T_1$-invariant subgroup. Set $T_2=\hat T_1|_{G_2}$. The group
$G_2$ acts by translations on $\hat G_1$
($g_2\cdot\chi_1=\chi_1+l(g_2)$), and we fall into the situation
described above (with $X=\hat G_1$, $G=G_2$, $T=\hat T_1$ and
$\beta=T_2$).

The roles of $G_1$ and $G_2$ above are almost symmetric. Indeed, to be
given an embedding $G_2\hookrightarrow\hat G_1$ with dense range is just the
same as to be given a non-degenerate pairing
$\langle\cdot\,,\,\cdot\rangle\colon G_1\times G_2\to\T$, then the equality
$T_2=\hat T_1|_{G_2}$ means that this pairing is $T_1\times T_2$-invariant.
The pairing gives rise to an embedding $r\colon G_1\hookrightarrow\hat G_2$.
Then $G_1$ acts on $\hat G_2$ by translations $g_1\cdot\chi_2=\chi_2-r(g_1)$,
and
$L^\infty(\hat G_1)\rtimes G_2\cong G_1\ltimes L^\infty(\hat G_2)$. In fact,
both algebras are canonically isomorphic to the twisted group W$^*$-algebra
$W^*(G_1\times G_2,\omega)$, where $\omega$ is the bicharacter defined by
$$
\omega((g'_1,g'_2),(g''_1,g''_2))=\langle g''_1,g'_2\rangle.
$$
Then $\alpha_T$ is nothing else than the automorphism induced by the
$\omega$-preserving automorphism $T_1\times T_2$.

\medskip

Let $R=\Z[t,t^{-1}]$ be the ring of Laurent polynomials over $\Z$,
$f\in\Z[t]$, $f\ne 1$, a polynomial whose irreducible factors are
not cyclotomic, equivalently, $f$ has no roots of modulus 1. Fix
$n\in\{2,3,\ldots,\infty\}$. Set $G_1=R/(f^\sim)$ and
$G_2=\mathop{\oplus}\limits^n_{k=1}R/(f)$, where
$f^\sim(t)=f(t^{-1})$. Let $T_i$ be the automorphism of $G_i$ of
multiplication by $t$. Let $\chi$ be a character of $G_2$. Then
the mapping $R\ni f_1\mapsto f_1(\hat T_2)\chi\in\hat G_2$ defines
an equivariant homomorphism $G_1\to\hat G_2$. In other words, if
$\chi=(\chi_1,\ldots,\chi_n)\in\hat G_2\subset\hat R^n$, then the
pairing is given by $$ \langle
f_1,(g_1,\dots,g_n)\rangle=\prod^n_{k=1}\chi_k(f^\sim_1\cdot g_k),
$$ where $(f^\sim_1\cdot g_k)(t)=f_1(t^{-1})g_k(t)$. This pairing
is non-degenerate iff the orbit of $\chi$ under the action of
$\hat T_2$ generates a dense subgroup of $\hat G_2$. Since $T_2$
is aperiodic, the dual automorphism is ergodic. Hence the orbit is
dense for almost every choice of $\chi$.

Now let us estimate entropy. First, by Yuzvinskii's formula \cite{Y,LW},
$h(\hat T_1)=m(f)$, $h(\hat T_2)=n\cdot m(f)$, where $m(f)$ is the
logarithmic Mahler measure of $f$,
$$
m(f)=\int^1_0\log|f(e^{2\pi is})|ds=\log|a_m|
   +\sum_{j\colon|\lambda_j|>1}\log|\lambda_j|,
$$
where $a_m$ is the leading coefficient of $f$, and $\{\lambda_j\}_j$ are
the roots of $f$. Now suppose that the coefficients of the leading
and the lowest terms of $f$ are equal to 1. Then $G_1\times G_2$ is a
free abelian group of rank $(n+1)\deg f$, and by a result of
Voiculescu~\cite{V} we have
$H(\alpha_T)\le h(\hat T_1\times\hat T_2)=(n+1)m(f)$.

Note also that since the automorphism $T_1\times T_2$ is aperiodic,
the automorphism $\alpha_T$ is mixing.

Let us summarize what we have proved.

\begin{theorem} \label{4.1}
For given $n\in\{2,3,\ldots,\infty\}$ and a polynomial $f\in\Z[t]$, $f\ne1$,
whose
coefficients of the leading and the lowest terms are equal to 1 and which
has no roots of modulus 1, there exists a mixing automorphism $\alpha$
of the hyperfinite II$_1$-factor and an $\alpha$-invariant Cartan subalgebra
$A$ such that
$$
H(\alpha|_A)=m(f),\ \ n\cdot m(f)\le H(\alpha)\le(n+1)m(f).
$$
\end{theorem}

The possibility of constructing on this way systems with arbitrary
values $H(\alpha|_A)<H(\alpha)$ is closely related to the
question, whether 0 is a cluster point of the set
$\{m(f)\,|\,f\in\Z[t]\}$ (note that it suffices to consider
irreducible polynomials whose leading coefficients and constant
terms are equal to 1). This question is known as Lehmer's problem,
and there is an evidence that the answer is {\it negative} (see
\cite{LSW} for a discussion).

\medskip

In estimating the entropy above we used the result of Voiculescu
stating that the entropy of an automorphism of a non-commutative
torus is not greater than the entropy of its abelian counterpart.
It is clear that this result should be true for a wider class of
systems. Consider the most simple case where the polynomial $f$ is
a constant.

\begin{example} \label{4.2}
\rm
Let $f=2$ and $n=2$. Then
$G_1=R/(2)\cong\mathop\oplus\limits_{k\in\Z}\Z/2\Z$, $G_2=G_1\oplus G_1$,
$T_1$ is the shift to the right, $T_2=T_1\oplus T_1$. Let
$G_1(0)=\Z/2\Z\subset G_1$ and $G_2(0)=\Z/2\Z\oplus\Z/2\Z\subset G_2$ be the
subgroups sitting at the 0th place. Set
$$
G^{(n)}_i=G_i(0)\oplus T_iG_i(0)\oplus\ldots\oplus T^n_iG_i(0).
$$
Then $\displaystyle
H(\alpha_T)\le hcpa_\tau(\alpha_T)\le\lim_{n\to\infty}{1\over n}\log\rank
   C^*(G^{(n)}_1\times G^{(n)}_2,\omega)\le 3\log2$, so (for
$A=L^\infty(\hat G_1)$)
$$
H(\alpha_T|_A)=\log2\ \ \hbox{and}\ \ 2\log2\le H(\alpha_T)\le 3\log2.
$$
The actual value of $H(\alpha_T)$ is probably depends on the choice of the
character $\chi\in\hat G_2$. We want to show that $H(\alpha_T)=2\log2$ for
some special choice of $\chi$. For this it suffices to require the
pairing
$\langle\cdot\,,\,\cdot\rangle|_{G^{(n)}_1\times G^{(n)}_2}$
be non-degenerate
in the first variable for any $n\ge0$ (so that $C^*(G^{(n)}_2)$ is a maximal
abelian subalgebra of $C^*(G^{(n)}_1\times G^{(n)}_2,\omega)$, and the rank of
the latter algebra is equal to $4(n+1)$).
The embedding $G_1\hookrightarrow\hat G_2$ is given by
$$
g_1\mapsto\prod_{n\in\Z:g_1(n)\ne0}\hat T^n_2\chi, \ \
 g_1=(g_1(n))_n\in\mathop{\oplus}_{n\in\Z}\Z/2\Z.
$$
So we must choose $\chi$ in a way such that the character
$\prod^m_{k=1}\hat T^{n_k}_2\chi$ is non-trivial on $G^{(n)}_2$ for any
$0\le n_1<\ldots< n_m\le n$. Identify $\hat G_2$ with
$\prod_{n\in\Z}(\Z/2\Z\oplus\Z/2\Z)$. Then $\hat T_2$ is the shift to the
right, and we may take any $\chi=(\chi_n)_n$ such that

(i) $\chi_n=0$ for $n<0$, $\chi_0\ne0$;

(ii) the group generated by $\hat T^n_2\chi$ is dense in $\hat G_2$.
\end{example}

Finally, we will show that it is possible to construct systems with positive
entropy, which have zero entropy on a Cartan subalgebra.

\begin{example} \label{4.3}
\rm Let $p$ be a prime number, $p\ne2$, $\hat G_1=\Z_p$ (the group
of $p$-adic integers), $G_2=\cup_{n\in\N}2^{-n}\Z\subset\hat G_1$,
$\hat T_1$ and $T_2$ the automorphisms of multiplication by $2$.
The group $G_1$ is the inductive limit of the groups $\Z/p^n\Z$,
and $T_1$ acts on them as the automorphism of division by~$2$.
Hence $$
H(\alpha_T|_A)=\lim_{n\to\infty}H(\alpha_T|_{C^*(\Z/p^n\Z)})=0. $$
Since $G_2=R/(t-2)$, we have $h(\hat T_2)=\log2$, so
$H(\alpha_T)\ge\log2$. We state that $$
H(\alpha_T)=hcpa_\tau(\alpha_T)=\log2. $$ The automorphism
$T^{p^{n-1}(p-1)}_1$ is identical on $\Z/p^n\Z$. Since $$
W^*(\Z/p^n\Z\times G_2,\omega)=\Z/p^n\Z\ltimes L^\infty(\hat G_2),
$$ by Theorem \ref{2.2} we infer $$
hcpa_\tau(\alpha^{p^{n-1}(p-1)}_T|_{W^*(\Z/p^n\Z\times
G_2,\omega)})
 =h(\hat T^{p^{n-1}(p-1)}_2),
$$
whence $hcpa_\tau(\alpha_T|_{W^*(\Z/p^n\Z\times G_2,\omega)})=\log2$, and
$$
hcpa_\tau(\alpha_T)=\lim_{n\to\infty}
   hcpa_\tau(\alpha_T|_{W^*(\Z/p^n\Z\times G_2,\omega)})=\log2.
$$
\end{example}

\bigskip\bigskip

\section{Flows on II$_1$-factors with invariant Cartan subalgebras} \label{5}

Using examples of previous sections and the construction of
associated flow we will construct systems with arbitrary values of
$H(\alpha|_A)$ and $H(\alpha)$ ($0\le H(\alpha|_A)\le H(\alpha)\le\infty$).

\medskip

Suppose a discrete amenable group $G$  acts freely and ergodically
by measure-preserving transformations $S_g$ on $(X,\mu)$, $T$ an
automorphism of $(X,\mu)$ and $\beta$ an automorphism of $G$ such
that $TS_g=S_{\beta(g)}T$. Consider the flow $F_t$ associated with
$T$. So $Y=\R/\Z\times X$, $d\nu=dt\times d\mu$, $$
F_t(\dot{r},x)=(\dot{r}+\dot{t},T^{[r+t]}x)\ \ \hbox{for}\ \
r\in[0,1),\
 x\in X,
$$
where $t\mapsto\dot{t}$ is the factorization mapping $\R\to\R/\Z$. The
semidirect product group
$G_0=G\times_\beta\Z$ acts on $(X,\mu)$. This action is ergodic. It is also
free, if
\begin{equation} \label{e5.1}
\hbox{there exist no }g\in G\hbox{ and no }n\in\N\hbox{ such that }S_g=T^n
 \hbox{ on a set of positive measure.}
\end{equation}
Let $\Gamma$ be a countable dense subgroup of $\R/\Z$, it acts by
translations on $\R/\Z$. Set $\G=\Gamma\times G_0$. The group $\G$
is amenable. It acts freely and ergodically on $(Y,\nu)$. The
corresponding equivalence relation is invariant under the flow, so
we obtain a flow $\alpha_t$ on $L^\infty(Y,\nu)\rtimes\G$. Compute
its entropy.
Let $\alpha_T$ be the automorphism of $L^\infty(X,\mu)\rtimes G$
defined by $T$. We state that
\begin{equation} \label{e5.3}
H(\alpha_t)=|t|H(\alpha_T),\ \ hcpa_\tau(\alpha_t)=|t|hcpa_\tau(\alpha_T),
   \ \ \hbox{and}\ \ H(\alpha_t|_{L^\infty(Y,\nu)})=|t|h(T).
\end{equation}
Since $h(F_t)=|t|h(F_1)=|t|h(\id\times T)$, the last equality in~(\ref{e5.3})
is evident. To prove the first two note that
$$
H(\alpha_t)=|t|H(\alpha_1)\ \ \hbox{and}\ \
 hcpa_\tau(\alpha_t)=|t|hcpa_\tau(\alpha_1)
$$
(see \cite[Proposition 10.16]{OP} for the first equality, the
second is proved analogously). We have
$$
L^\infty(Y,\nu)\rtimes\G=(L^\infty(\R/\Z)\rtimes\Gamma)\otimes
   (L^\infty(X,\mu)\rtimes G_0),\ \ \alpha_1=\id\otimes\tilde\alpha_T,
$$ where $\tilde\alpha_T$ is the automorphism of
$L^\infty(X,\mu)\rtimes G_0$ defined by $T$. Since completely
positive approximation entropy is subadditive and monotone
\cite{V}, we have
$hcpa_\tau(\id\otimes\tilde\alpha_T)=hcpa_\tau(\tilde\alpha_T)$.
We have also $H(\id\otimes\tilde\alpha_T)=H(\tilde\alpha_T)$ by
Remark following the proof of Theorem~\ref{2.2}. Since $$
L^\infty(X,\mu)\rtimes
G_0=(L^\infty(X,\mu)\rtimes_SG)\rtimes_{\alpha_T}\Z, $$ we obtain
$hcpa_\tau(\tilde\alpha_T)=hcpa_\tau(\alpha_T)$ and
$H(\tilde\alpha_T)=H(\alpha_T)$ by virtue of Theorem~\ref{2.2}. So
$hcpa_\tau(\alpha_1)=hcpa_\tau(\alpha_T)$ and
$H(\alpha_1)=H(\alpha_T)$,
and the proof of the equalities~(\ref{e5.3}) is complete.

\begin{theorem} \label{5.1}
For any $s$ and $t$, $0\le s< t\le\infty$, there exist an automorphism
$\alpha$ of the hyperfinite II$_1$-factor and an $\alpha$-invariant Cartan
subalgebra~$A$ such that
$$
H(\alpha|_A)=s\ \ \hbox{and} \ \ H(\alpha)=t.
$$
\end{theorem}

\begin{verif}{}{}
Consider a system from Example~\ref{4.3}. Then the condition
(\ref{e5.1}) is satisfied, so the construction above leads to a
flow $\alpha_t$ and an $\alpha_t$-invariant Cartan subalgebra
$A_1$ such that $$ H(\alpha_t|_{A_1})=0 \ \ \hbox{and}\ \
H(\alpha_t)=hcpa_\tau(\alpha_t)
 =|t|\log2.
$$
As in Example~\ref{3.3}, consider a Bernoulli flow $S_t$ on $(X,\mu)$ with
$h(S_1)=\log2$. Then for the corresponding flow $\beta_t$ on
$L^\infty(X,\mu)\rtimes_{S_1}\Z$ we have (with $A_2=L^\infty(X,\mu)$)
$$
H(\beta_t|_{A_2})=H(\beta_t)=hcpa_\tau(\beta_t)=|t|\log2.
$$
Since Connes-St{\o}rmer' entropy is superadditive~\cite{SV} and Voiculescu's
entropies are subadditive, we conclude that
$$
H((\alpha_t\otimes\beta_s)|_{A_1\otimes A_2})=|s|\log2,\ \
 H(\alpha_t\otimes\beta_s)=H(\alpha_t)+H(\beta_s)=(|t|+|s|)\log2.
$$

Finally, consider an infinite tensor product of systems from
Example~\ref{4.3}. Thus we obtain an automorphism $\gamma$ and an
$\alpha$-invariant Cartan subalgebra~$A_3$ such that
$$
H(\gamma|_{A_3})=0\ \ \hbox{and}\ \ H(\gamma)=\infty.
$$
Then $H(\beta_s\otimes\gamma)|_{A_2\otimes A_3})=|s|\log2$,
$H(\beta_s\otimes\gamma)=\infty$.
\end{verif}

\bigskip\bigskip

\section{Final remarks} \label{6}

\subsection
Let $p_1$ and $p_2$ be prime numbers, $p_i\ge3$,
$i=1,2$. Construct automorphisms $\alpha_1$ and $\alpha_2$ as in
Example~\ref{4.3}.

\begin{proposition} If $p_1\ne p_2$, then $\alpha_1$ and $\alpha_2$ are not
conjugate as automorphisms of the hyperfinite II$_1$-factor,
though $H(\alpha_1)=H(\alpha_2)=\log 2$.
\end{proposition}

\begin{verif}{}
Indeed, the automorphisms define
unitary operators $U_i$ on $L^2(M, \tau)$.
As we can see, the point part $S_i$ of the
spectrum of $U_i$ is non-trivial. If $p_1\ne p_2$, then
$S_1\ne S_2$, so $\alpha_1$ and $\alpha_2$ are not
conjugate.
\end{verif}

\subsection
The automorphisms of Theorem~\ref{4.1} and Example~\ref{4.2} are
ergodic. On the other hand, the automorphisms of Example~\ref{4.3}
are not ergodic, even on the Cartan subalgebra. Moreover, any
ergodic automorphism of compact abelian group has positive entropy
(it is even Bernoullian), so with the methods of Section~\ref{4}
we can not construct ergodic automorphisms with positive entropy
and zero entropy restriction to a Cartan subalgebra (however, for
actions of $\Z^d$, $d\ge2$, we are able to construct such
examples).

The construction of Section~\ref{5} leads to non-ergodic automorphisms also,
even if we start with an ergodic automorphism (such as in Example~\ref{4.2}).

\medskip

{\bf Acknowledgement.} The first author (V.G.) is grateful to Erling
St{\o}rmer for interesting and helpful discussions of the first version
of this paper.

\end{document}